# Fast and Accurate Optimization of Metasurfaces with Gradient Descent and the Woodbury Matrix Identity

Jordan Budhu, *Member, IEEE*, and Anthony Grbic, *Fellow, IEEE*

*Abstract*—A fast metasurface optimization strategy for finite-size metasurfaces modeled using integral equations is presented. The metasurfaces considered are constructed from finite patterned metallic claddings supported by grounded dielectric spacers. Integral equations are used to model the response of the metasurface to a known excitation and solved by Method of Moments. An accelerated gradient descent optimization algorithm is presented that enables the direct optimization of such metasurfaces. The gradient is normally calculated by solving the method of moments problem $N + 1$ times where $N$ is the number of homogenized elements in the metasurface. Since the calculation of each component of the $N$-dimensional gradient involves perturbing the moment method impedance matrix along one element of its diagonal and inverting the result, this numerical gradient calculation can be accelerated using the Woodbury Matrix Identity. The Woodbury Matrix Identity allows the inverse of the perturbed impedance matrix to be computed at a low cost by forming a rank-r correction to the inverse of the unperturbed impedance matrix. Timing diagrams show up to a 26.5 times improvement in algorithm times when the acceleration technique is applied. An example of a passive and lossless wide-angle reflecting metasurface designed using the accelerated optimization technique is reported.

## I. INTRODUCTION

METASURFACES are subwavelength textured surfaces which perform desired wavefront transformations [1]. They can be designed using a boundary condition called the Generalized Sheet Transition Condition or GSTC [2,3] and realized as a cascade of sheets and dielectrics [4-8]. The cascade is designed to generate the same transmission and reflection properties as the GSTC boundary. For homogeneous structures under normal incidence, the approach works well and has resulted in numerous designs involving polarization rotation [9,10] and polarization conversion [11]. More complex field transformations are enabled by incorporating spatial inhomogeneity such as wide-angle refraction [7,12] or beam collimation [13,14]. The same approaches have been adopted to design these inhomogeneous metasurfaces. However, these techniques can lead to inaccuracies since these simplified transmission-line and network-based theories do not model transverse coupling between the cells of an inhomogeneous metasurface. To minimize transverse coupling, the metasurfaces had to be restrictively thin [15,16] or contain conducting baffles separating the unit cells [17,18]. A new approach to design more practical metasurfaces involving integral equations solves the problem [19-21]. Rather than using the infinitely-thin GSTC boundary condition in design, the actual thickness of the metasurface is modeled in addition to the transvers coupling between elements. Since the integral equation approach models the entire metasurface rather than a single unit cell, optimization of all the elements of a finite inhomogeneous metasurface becomes possible [20,21]. Approaches which model the entire metasurface geometry have been reported by others. In [22,23], the Rigorous Coupled Wave Analysis (RCWA) was applied in conjunction with adjoint variable optimization methods to optimize very large metasurfaces. Although the technique optimizes an inhomogeneous metasurface, the metasurface was broken down into cells a few wavelengths across and local periodicity applied

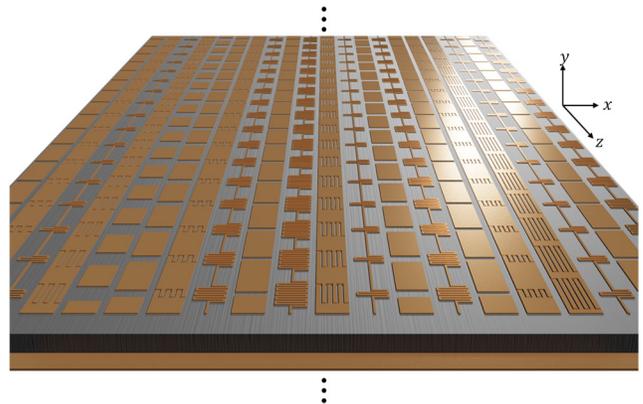

Figure 1. Metasurface geometry. The metasurface is infinite in length and invariant in the $z$ direction and finite in width and spatially variant in the $x$ direction. The metasurface consists of three layers: one patterned metallic cladding layer separated from a ground plane by a dielectric spacer.

This work was supported in part by the Office of Naval Research under grant no. N00014-18-1-2536 and the Army Research Office under grant no. W911NF-19-1-0359.

J. Budhu and A. Grbic are with the is with the Department of Electrical Engineering and Computer Science at the University of Michigan, Ann Arbor, MI (e-mail: jbudhu@umich.edu, agrbic@umich.edu).
.

to the analysis of those cells. Using this approach, they were able to optimize multilayered metasurfaces over $1000\lambda$ in diameter. This approach, however, does not model the finite metasurface width and will be approximate at interfaces where the geometry in two neighboring unit cells vary greatly. In [24], gradient descent optimization was accelerated with the adjoint variable method to optimize large metasurfaces and compared to global optimizers such as Genetic Algorithm. The adjoint variable method is very powerful and can calculate gradients rapidly. The authors report optimization of metasurfaces a few tens of wavelengths in size using this approach. However, when the metasurfaces are larger, the authors preoptimize part libraries which can be stitched together to form the optimized final metasurface. In [25], the Alternating Direction Method of Multipliers technique is applied to optimize bianisotropic metasurfaces modeled using integral equations. Their approach, however, does not model the finite thickness of the metasurface or the dielectric spacers which would support them. In [26], gradient descent optimization was applied to obtain reactive mode converting metasurfaces placed within a cylindrical waveguide. However, this approach is geometry specific, and therefore is restricted to modes with radial variations.

This paper presents an optimization approach which models and optimizes an entire finite-sized metasurface. With the possibility of optimizing large inhomogeneous metasurfaces with finite dimensions and thicknesses, new designs which perform more complex field transformations are possible [27]. It is based on a gradient descent optimization technique accelerated by the Woodbury Matrix Identity [28]. The Woodbury Matrix Identity (also known as the Sherman-Morrison-Woodbury formula or the matrix inversion lemma) says that the inverse of a rank-r correction of some matrix can be computed by doing a rank-r correction to the inverse of the original matrix [28]. When metasurfaces are modelled using integral equations, the cost function is evaluated by solving a matrix equation. Each component of the gradient of the cost function involves forming a rank-r correction to the moment method impedance matrix. This is the case since only one diagonal element of the impedance matrix is perturbed per dimension. Thus, the perturbed matrix inverse can be calculated by forming a rank-r correction to the unperturbed matrix inverse. This identity is commonly used to accelerate method of moments problems in which small perturbations are made to an impedance matrix [29,30]. The identity has also been applied to the simulation of power electronics circuits which can be formulated in terms of a matrix equation [31].

This paper is organized as follows. Section II presents the accelerated gradient descent optimization algorithm. The Woodbury Matrix Identity is presented along with the formulations needed to optimize metasurfaces. In section III, timing diagrams and data are provided comparing the algorithm time for both the direct gradient calculation and the calculation accelerated with the Woodbury Matrix Identity. In section IV, an example of a finite-sized, wide-angle reflecting metasurface is provided. Finally, in section V, the paper is concluded.

## II. METASURFACE DESIGN AND OPTIMIZATION

The metasurfaces considered consist of a patterned metallic cladding (layer 1) supported by a dielectric spacer (layer 2) and backed by a conducting ground plane (layer 3), as shown in Fig. 1. The metasurface and supporting grounded dielectric substrate are finite-sized in the $x$ and $y$ directions, and infinite in the $z$ direction. Therefore, the geometry is two-dimensional. The metasurfaces are designed following the integral equation design technique outlined in [20]. To summarize, the metasurface elements are homogenized (see Fig. 2) and the metasurface is modeled as an inhomogeneous array of homogeneous sheet impedances. The response of the metasurface to an excitation is modeled using integral equations. The sheet impedances are computed by directly solving the governing integral equations by the method of moments. The method of moments solution yields complex-valued sheet impedances with possibly both positive and negative resistances if the local power density of the incident field differs from that of the scattered field at the metasurface [20,21,32-34]. Since a passive and lossless metasurface is desired, the sheet impedances are optimized to obtain purely reactive responses, thereby introducing a number of surface waves which carry power across the metasurface.

In gradient-based optimization methods, convergence strongly depends on obtaining a good initial solution. A good initial solution can be obtained from the complex-valued impedance sheet obtained by directly solving the governing integral equation through the Method of Moments. The real parts of the complex-valued sheet impedances are discarded and the reactances kept. These $N$ unknown reactances are arranged in an $N$-dimensional space with each reactance varying along an orthogonal axis. Since only reactances are retained in the optimization, the final design is passive and lossless. A surface is defined in this space as $f(\vec{x}) = f(x_1, x_2, x_3, \ldots, x_N)$ and represents the response of the metasurface as a function of its reactances. The cost function $f$ is designed such that its minimum represents the optimal solution.

### A. Newton Optimization Method

At an initial point $\vec{x}_t$, the function $f$ is expanded into a second order Taylor series expansion

$$f(\vec{x}) \approx f(\vec{x}_t) + \left[\nabla f(\vec{x}_t)\right]^T (\vec{x} - \vec{x}_t) + \frac{1}{2}(\vec{x} - \vec{x}_t)^T H(\vec{x}_t)(\vec{x} - \vec{x}_t) \quad (1)$$

where $\nabla f(\vec{x})$ is the gradient vector and $H(\vec{x}_t)$ is the Hessian matrix of second partial derivatives of $f$. When $H$ is symmetric and positive definite, the quadratic approximation has a well-defined minimum at

$$\nabla f(\vec{x}) \approx \nabla f(\vec{x}_t) + H(\vec{x}_t)(\vec{x} - \vec{x}_t) = 0 \quad (2)$$

obtained by setting the gradient of (1) equal to zero. Solving (2) for $\vec{x}$ gives the next iterate $\vec{x}_{t+1}$

$$\vec{x}_{t+1} = \vec{x} - H^{-1}(\vec{x}_t) \nabla f(\vec{x}_t) \quad (3)$$

In practice, (3) generally includes a controllable parameter $\alpha$ that limits the size of the step

$$\vec{x}_{t+1} = \vec{x} - \alpha H^{-1}(\vec{x}_t) \nabla f(\vec{x}_t) \quad (4)$$

where $0 < \alpha$. The function $f$ is expanded about the new point $\vec{x}_{t+1}$ in a second order Taylor series using (1) again and the process is repeated to obtain the next iterate $\vec{x}_{t+2}$ and so forth until the minimum is reached or the algorithm converges (see Section II.E for convergence criteria). Thus, at each point $\vec{x}$, the gradient of $f$, the Hessian of $f$, and the step size $\alpha$ must be calculated.

*B. Calculation of Hessian*

The Hessian can be approximated from successive calculations of the gradient $\nabla f(\vec{x})$ using the secant method

$$H(\vec{x}_t) = \nabla^2 f(\vec{x}_t) \approx \frac{\nabla f(\vec{x}_t) - \nabla f(\vec{x}_{t-1})}{\vec{x}_t - \vec{x}_{t-1}} \quad (5)$$

The inverse of $H$ for the next iterate, $H_{t+1}^{-1}$, can be formed from the inverse of $H$ at the current iterate, $H_t^{-1}$, from [35]

$$H_{t+1}^{-1} = \left(I_{kN} - \frac{sy^T}{y^Ts}\right) H_t^{-1} \left(I_{kN} - \frac{ys^T}{y^Ts}\right) + \frac{ss^T}{y^Ts} \quad (6)$$

where $I_{kN}$ is the $kN \times kN$ identity matrix, $s = \vec{x}_t - \vec{x}_{t-1}$, and $y = \nabla f(\vec{x}_t) - \nabla f(\vec{x}_{t-1})$. The Hessian in (5) and its inverse calculated using (6) is guaranteed to be symmetric and positive definite. For the first iteration, the Hessian can be initialized to $I_{kN}$. Approximating the Hessian inverse in this way in (6) is referred to as a Quasi-Newton Algorithm. Thus, all of the computational expense in calculating the next iterate using (4) lies in the calculation of the gradient vector.

*C. Calculation of the Numerical Gradient*

The numerical gradient is the most costly part of (4) to calculate. The gradient is approximated using forward finite differences

$$\nabla f(\vec{x}_t) = \begin{bmatrix} \frac{\partial f(\vec{x}_t)}{\partial x_1} \\ \frac{\partial f(\vec{x}_t)}{\partial x_2} \\ \vdots \\ \frac{\partial f(\vec{x}_t)}{\partial x_N} \end{bmatrix} = \begin{bmatrix} \frac{f(x_1+\delta_1, x_2,...,x_N) - f(\vec{x}_t)}{\delta_1} \\ \frac{f(x_1, x_2+\delta_2,...,x_N) - f(\vec{x}_t)}{\delta_2} \\ \vdots \\ \frac{f(x_1, x_2,...,x_N+\delta_N) - f(\vec{x}_t)}{\delta_N} \end{bmatrix} \quad (7)$$

To evaluate the gradient in (7), $N+1$ evaluations of $f$ are required. If $f$ is computationally expensive, then the gradient vector may be the temporal bottleneck of the optimization algorithm.

When the metasurface response is modeled using integral equations [20] and discretized, the evaluation of $f$ involves the solution of a matrix equation

$$[I] = [Z]^{-1}[V] \quad (8)$$

where $[I]$, $[Z]$, and $[V]$ are the method of moment current, impedance, and voltage vectors, respectively. For example, the far-field radiation pattern resulting from the currents in (8) can be calculated and the cost function $f$ can be the root-mean-square (RMS) difference between this pattern and a target pattern (see (19) for example). Since the metasurface discretization is usually subwavelength, to adhere to the homogenization limit, the size of the matrix $[Z]$ can be large and the calculation of its inverse time consuming. The first step in the calculation of the gradient in (7) is to calculate $f(\vec{x}_t)$ by forming the impedance matrix as

$$[Z] = [Z_c] + [\eta_s] \quad (9)$$

and solving (8). In (9), $[Z_c]$ is the mutual coupling matrix between the homogenized elements of the metasurface, and $[\eta_s]$ is a diagonal matrix with the sheet reactances $\vec{x}_t$ defined along the diagonal. For example, in [20, Eqn. (22)], the following matrix equation is defined

$$\begin{bmatrix} [V_m^M] \\ [V_m^G] \end{bmatrix} = \begin{bmatrix} [Z_{mn}^{MM}] + [\eta_s] & [Z_{mn}^{MG}] \\ [Z_{mn}^{GM}] & [Z_{mn}^{GG}] \end{bmatrix} \begin{bmatrix} [I_n^M] \\ [I_n^G] \end{bmatrix} \quad (10)$$

where the matrices $[Z_c]$ and $[\eta_s]$ are

$$[Z_c] = \begin{bmatrix} [Z_{mn}^{MM}] & [Z_{mn}^{MG}] \\ [Z_{mn}^{GM}] & [Z_{mn}^{GG}] \end{bmatrix}, \quad [\eta_s] = \begin{bmatrix} [\eta_s] & [0] \\ [0] & [0] \end{bmatrix} \quad (11)$$

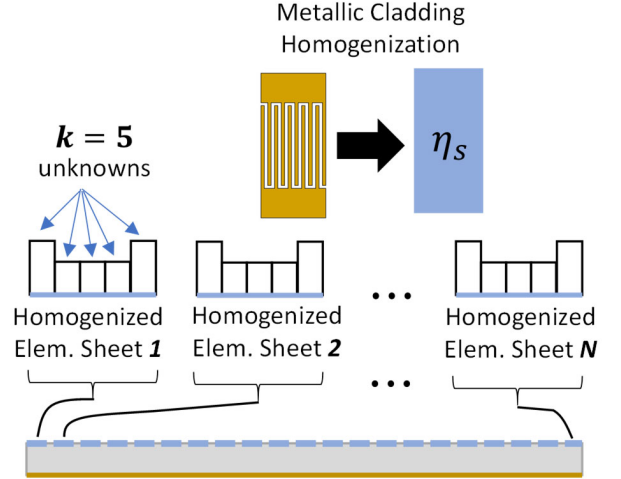

Figure 2. Metasurface homogenization and placement of current unknowns. Each of the $N$ elements of the patterned metallic cladding layer of the metasurface are homogenized and assigned a single sheet impedance. $k$ pulse basis functions are then placed along each of the $N$ homogenized sheet impedances. The total number of unknowns is then $kN$.

The matrix equation (10) models an impedance sheet representing a metasurface (superscript M) above a perfectly conducting ground plane (superscript G). The matrices $[V_m^{M,G}]$ are therefore the voltage vectors corresponding to the metasurface and the ground plane, respectively. Similarly, the impedance matrix $[Z_{mn}^{MG}]$ models the mutual coupling between currents on the metasurface and the ground plane.

The next step in the calculation of each element of the gradient in (7) involves a small perturbation of one reactance sheet ($f(x_1 + \delta_1, x_2, ... x_N)$ for example). To evaluate $f$ for the perturbed vector, the matrix equation (8) must be solved again using the perturbed vector in (9). The matrix appearing in (8) is of dimension $kN \times kN$, where $N$ is the number of sheet elements in the metasurface and $k$ is the number of unknowns placed on each of those $N$ elements (see Fig. 2).

The Woodbury matrix identity can now be used to form the matrix inverse of the perturbed $kN \times kN$ matrix $[Z]$ appearing in (8) without having to calculate the $kN \times kN$ inverse. The Woodbury matrix identity states that a low rank update to a matrix inverse can be found using

$$(Z + ACB)^{-1} = Z^{-1} - Z^{-1}A(C^{-1} + BZ^{-1}A)^{-1} BZ^{-1} \quad (12)$$

where $Z$ is the $kN \times kN$ $[Z]$ matrix defined in (9), $C$ is a $k \times k$ matrix which contains the rank $r$ update to the matrix $Z$, $A$ and $B$ are a $kN \times k$ matrix and a $k \times kN$ matrix used to form the $kN \times kN$ rank $r$ matrix update to matrix $Z$. Equation (12) thus permits the evaluation of the perturbed $f$ inherent in each element of the gradient in (7) by inverting only a $k \times k$ matrix $(C^{-1} + BZ^{-1}A)$ and forming the matrix products of (12). Since the matrix $C$ is diagonal, its inverse is formed analytically by inverting its diagonal elements. When the metasurface contains a large number of elements, $N$ can be hundreds or thousands, and $k$ may be up to ten or more depending on how many unknowns are placed per reactive element in the metasurface (See Fig. 2). Since calculation of the matrix inverse is the costliest part of evaluation of $f$, the Woodbury matrix identity can accelerate the gradient calculation by orders of magnitude.

The matrices $A$ and $B$ are rectangular matrices which determine the position of the elements of $C$ within the $kN \times kN$ product matrix $ACB$. For example, if the second element of the gradient vector is being calculated and there are $k = 5$ unknowns per metasurface element (See Fig. 2), then the matrices $A, C,$ and $B$ would be defined as

$$A = \begin{bmatrix} [0]_{5 \times 5} \\ [I]_{5 \times 5} \\ [0]_{(5N-10) \times 5} \end{bmatrix}, \quad B = A^T$$

$$C = \begin{bmatrix} \delta_2 & 0 & 0 & 0 & 0 \\ 0 & \delta_2 & 0 & 0 & 0 \\ 0 & 0 & \delta_2 & 0 & 0 \\ 0 & 0 & 0 & \delta_2 & 0 \\ 0 & 0 & 0 & 0 & \delta_2 \end{bmatrix} \quad (13)$$

where $[I]$ is the identity matrix. In general, if $k$ unknowns are placed on each of the $N$ elements of the metasurface (See Fig. 2), the following matrices correspond to the $l^{th}$ sheet element

$$A = \begin{bmatrix} [0]_{k(l-1) \times k} \\ [I]_{k \times k} \\ [0]_{(kN-kl) \times k} \end{bmatrix}, \quad B = A^T, \quad C = \delta_l [I]_{k \times k} \quad (14)$$

Since the matrices $A$ and $B$ are sparse, fast matrix multiplication routines for the matrix products in (12) can further improve the algorithm time. Using (12) and (14), the gradient vector (7) can be calculated rapidly. In each iteration, the matrix inverse calculations needed to compute the gradient vector can be found from low-cost computations involving the Woodbury matrix identity. Hence, finite-sized non-periodic metasurfaces with large numbers of unknowns can be optimized in feasible times which may not be possible without these techniques.

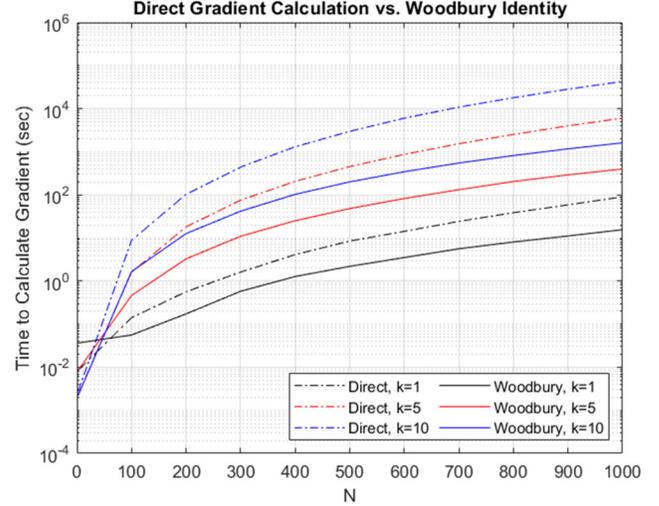

Figure 3. Direct gradient calculation versus Woodbury Matrix Identity timing plot. The horizontal axis represents the number of unknowns and the vertical axis the time required to calculate the full gradient vector. The different curves are for different $k$ values (See Fig. 2). The dashed curves are for the direct calculation method and the solid curves are for the Woodbury Matrix Identity method.

### D. Calculation of Step Size

Once the Hessian and gradient are calculated from (6) and (7), the step size is calculated by performing a line search along the Newton direction $d = -H^{-1}(\vec{x}_t)\nabla f(\vec{x}_t)$. The line search is a 1-dimensional optimization

$$\alpha = \arg \min_{\alpha \geq 0} f(\vec{x}_t + \alpha d) \quad (15)$$

A golden section optimization is used to perform the 1-dimensional line search optimization in (15) [36].

### E. Convergence Criteria

The optimization is ended when either one of three criteria are met.

1) the norm of the gradient vector, $|\nabla f|$, falls below some threshold, $a_1$

$$|\nabla f| < a_1$$

2) the step size, $\alpha$, falls below a threshold $a_2$

$$\alpha < a_2$$

3) the evaluation of $f$ falls below a threshold $a_3$

$$f < a_3$$

## III. ACCELERATED GRADIENT CALCULATION RESULTS

In this section, the gradient (7) is calculated both directly, i.e. without using the Woodbury Matrix Identity ('Direct'), and by using the Woodbury Matrix Identity ('Woodbury') for various combinations of $N$ and $k$. The results are shown in Table 1 and also plotted in Fig. 3. The data was obtained from the calculation of the inverse of a randomly generated impedance matrix with the dimensions indicated in the table.

From the data, we observe up to a 26.5 times improvement in time to calculate the gradient vector. Since the entire optimization algorithm rests on calculating the gradient vector, the optimization time is approximately represented by the presented data in Table 1 and Fig. 3. In the next section, we present a design example of wide-angle reflection metasurface designed using the presented accelerated optimization technique.

TABLE I
TIME TO CALCULATE GRADIENT TABLE

| N | k=10 Woodbury | k=10 Direct | k=5 Woodbury | k=5 Direct | k=1 Woodbury | k=1 Direct |
|---|---|---|---|---|---|---|
| 2 | 0.00219 | 0.00275 | 0.00851 | 0.00228 | 0.03655 | 0.00842 |
| 100 | 1.65428 | 8.75606 | 0.46304 | 1.60233 | 0.05596 | 0.14142 |
| 200 | 12.6309 | 102.928 | 3.27053 | 18.1083 | 0.17314 | 0.56446 |
| 300 | 41.8161 | 440.721 | 10.9370 | 75.2882 | 0.57928 | 1.61378 |
| 400 | 102.723 | 1303.60 | 25.0543 | 205.661 | 1.27434 | 4.11161 |
| 500 | 199.822 | 2976.71 | 48.0400 | 453.314 | 2.19397 | 8.44913 |
| 600 | 343.984 | 6021.33 | 82.2661 | 868.790 | 3.50960 | 14.2327 |
| 700 | 547.354 | 10861.4 | 131.575 | 1537.37 | 5.60646 | 24.2042 |
| 800 | 815.152 | 17912.7 | 203.541 | 2515.70 | 8.05453 | 38.5044 |
| 900 | 1167.28 | 28221.8 | 290.394 | 3995.00 | 11.1107 | 58.5054 |
| 1000 | 1610.26 | 42697.1 | 397.627 | 6027.14 | 15.5304 | 89.8736 |

*Note, all times are in seconds. The shaded columns are calculated with the Woodbury Matrix identity while the unshaded columns are not. The times in the table were calculated on a machine with a 6-core 2.2-GHz Intel Core i7-8750H CPU.

## IV. PERFECT WIDE-ANGLE REFLECTION METASURFACE EXAMPLE

This section presents the design of a wide-angle perfectly reflecting metasurface. The metasurface considered consists of an impedance sheet over a grounded dielectric substrate (see Fig. 4). The metasurface is $8\lambda_0$ wide at 10GHz. The Rogers 6010 substrate has permittivity $\epsilon_r=10.7(1-j0.0023)$ and thickness 0.00127m ($\lambda_0/23.62$). The impedance sheet has been divided into 160 cells of width $\lambda_0/20$. The metasurface is designed to reflect a normally incident plane wave to a reflection angle of 70°. Thus, the incident field is

$$\vec{E}^{inc} = E_0 e^{jk_0 y} \hat{z} \quad (16)$$

The desired scattered field is

$$\vec{E}^{sca} = E_0 e^{jk_0(x\sin 70° - y\cos 70°)} \hat{z} \quad (17)$$

where the reflected amplitude has been set equal to the incident field amplitude ($E_r = E_0$). This choice will lead to a perfectly reflecting metasurface which is lossy [38]. Note, alternatively, one could start with the global power conserving choice of $E_r = E_0 \sqrt{\cos(\theta_i)/\cos(\theta_r)}$, which will lead to a reflecting metasurface which exhibits equal amounts of loss and gain.

The metasurface is modelled using the integral equation technique outlined in [20]. The integral equations are transformed into linear matrix equations using the method of moments. The linear matrix equation for this case is

$$\sum_{i=1}^{3}\sum_{j=1}^{3}\left([V_i] = [Z_{ij}][I_j] + [\eta_i][I_i]\right) \quad (18)$$

where $i=1,2,3$ denotes the different observation layers and $j=1,2,3$ denotes the different source layers (see Fig. 4). Note, in (18), $[\eta_{s3}]=0$ since the ground plane is assumed perfectly conducting, and $[\eta_{s2}] = [j\omega\epsilon_0(\epsilon_r-1)]^{-1}$ for the substrate by the volume equivalence principle. Directly solving (18) for $[\eta_{s1}]$ yields the complex-valued sheet impedances shown in Fig. 5a. Since the local power density is not conserved at each point in layer 1, the sheets impedances are complex-valued. Had the metasurface been infinite in extent, the real part of Fig. 5a would be completely lossy (positive resistances) due to the choice of (17). However, since the metasurface is finite in the $x$-direction, some gain is required to cancel edge diffraction.

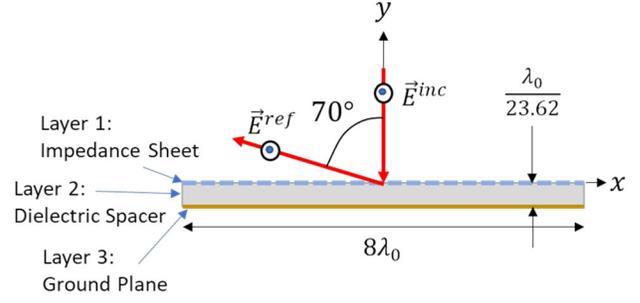

Figure 4. Metasurface geometry. The metasurface contains 3 layers (an impedance sheet, a dielectric spacer, and a ground plane). The geometry is invariant in the z-direction and thus the problem is 2-dimensional.

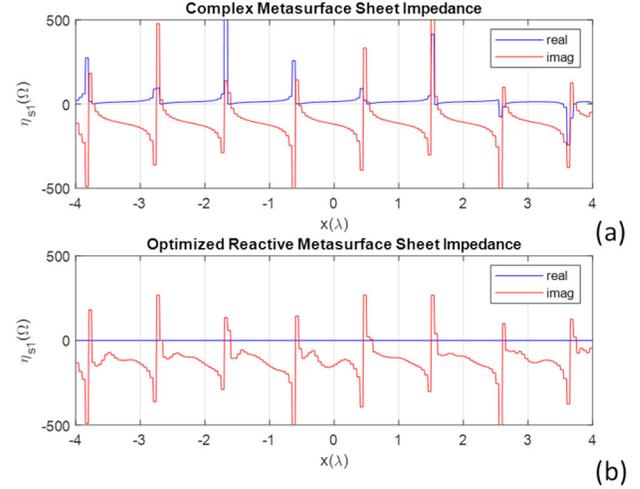

Figure 5. Metasurface layer 1 impedance sheet. (a) The complex-valued sheet resulting from the direct solution of the matrix equation. (b) The optimized purely reactive sheet.

Next, the complex $[\eta_{s1}]$ will be transformed into a purely reactive $[\eta_{s1}]$ using the accelerated gradient descent optimization from section II. The optimization domain is 160 dimensional with one reactance varying along each dimension. Ten unknowns were placed on each reactance. Therefore, $N=160$ and $k=10$. The initial point (seed solution) in the domain is formed from the imaginary parts of Fig. 5a. An optimization cost function $f$ is defined in this space as a function of the reactances of layer 1 as

$$f([\eta_{s1}]) = RMS\left\{\left|E_{farfield}(\phi)\right|_{\substack{\text{Complex}\\\text{Sheet}}} - \left|E_{farfield}(\phi)\right|_{\substack{\text{Reactive}\\\text{Sheet}}}\right\} \quad (19)$$

where $|E_{farfield}(\phi)|$ is the normalized magnitude of the scattered far field pattern. It is calculated by solving (18) for the induced currents $[I_1], [I_2],$ and $[I_3]$ (the surface current density along layer 1, the polarization current density within layer 2, and the surface current density along layer 3, respectively) then computing the far field pattern using the magnetic vector potential formulation. Note, RMS is the root-mean-square value of the quantity in the brackets. The magnitude of the scattered far field pattern resulting from the metasurface with the complex $[\eta_{s1}]$, $|E_{farfield}(\phi)|_{complex\,sheet}$, is shown in Fig. 6 as the curve labeled 'MoM Complex'. The cost function minimizes the difference between the normalized magnitude of the far field pattern radiated by the reactive sheet and the complex sheet.

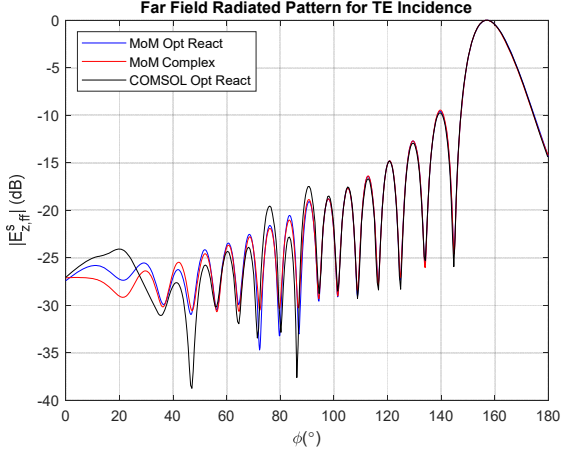

Figure 6. Scattered far field patterns. Shown is the far field pattern resulting from the complex-valued sheet calculated using the MoM (MoM Complex), the far field resulting from the optimized reactive sheet calculated using the MoM (MoM Opt React), and the far field pattern resulting from the optimized reactive sheet calculate in COMSOL Multiphysics (COMSOL Opt React).

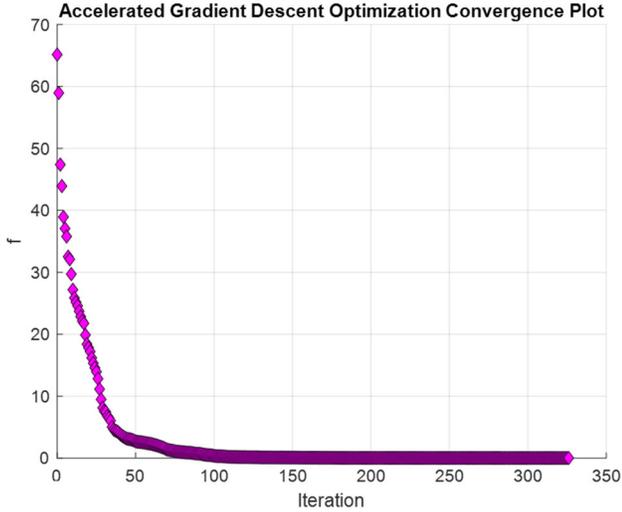

Figure 7. Optimization convergence plot.

The optimization converged in 326 iterations to a value of $f$=0.022 (see Fig. 7) because the step size $\alpha < 10^{-10}$ (see section II.E.2). The optimization took 4 hours and 51 minutes running on 25 total cores of a machine with dual 3.0 GHz Intel Xeon Gold 6154 CPUs. According to Table I and Fig. 3, for $N$=160 and $k$=10, the 'Woodbury' algorithm runs approximately 8 times faster than the 'Direct' algorithm. Thus, without the acceleration, the optimization would have taken 38 hours and 48 minutes (over 1.5 days).

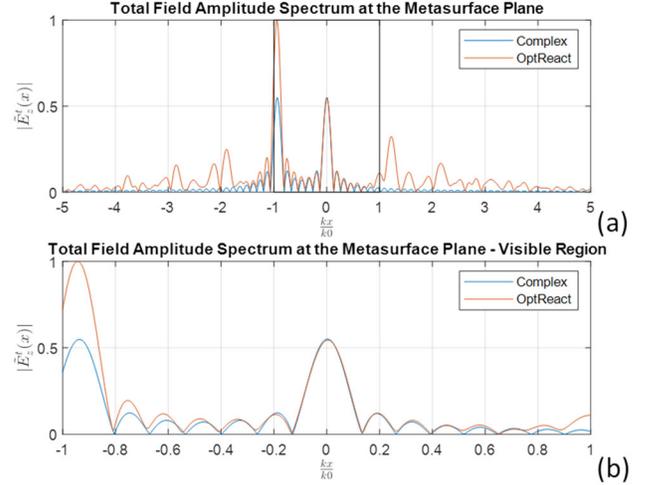

Figure 8. Total field amplitude spectra. (a) The visible and invisible spectrum. (b) Only the visible spectrum.

The optimized reactive $[\eta_{s1}]$ is shown in Fig. 5b. Note the impedance is purely reactive with no real part. Furthermore, the optimization has introduced perturbations to the imaginary part to excite evanescent waves which together with the incident and scattered propagating fields, satisfy local conservation of power density across the sheet. These evanescent waves are evident in Fig. 8a. The figure shows the amplitude spectrum of the total electric field on layer 1 for both the complex-valued metasurface and the optimized purely reactive metasurface. The plot is normalized to the scattered peak for the optimized purely reactive metasurface. As can be seen, the optimized reactive metasurface supports significant evanescent spectrum whereas the complex valued sheet does not. Note also that the result for the complex-valued metasurface shows equal peak amplitude in both spectral peaks at $k_x/k_0$=0 and $k_x/k_0$=-0.94 since $E_r = E_0$ in (17). The peak centered at $k_x/k_0$=0 is due to the incident field and the peak centered at $k_x/k_0$=-0.94 is due to the scattered propagating field. For the optimized purely reactive metasurface, the amplitude of the scattered field is greater than that of the incident field, in order to satisfy global conservation of power density in a passive way. For perfect anomalous reflection, the power carried in the desired reflection direction must be equal to the power of the incident plane wave. This requires the reflected plane wave to carry amplitude [38]

$$|E^{sca}| = \frac{\sqrt{\cos\theta_i}}{\sqrt{\cos\theta_r}} = 1.7 \text{ V/m} \quad (20)$$

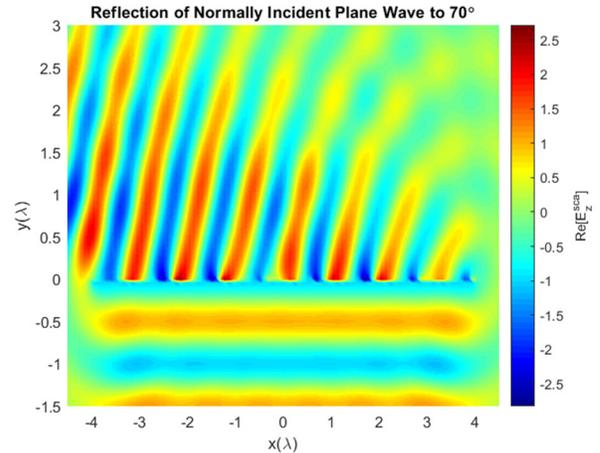

Figure 9. Real part of the complex scattered electric near field.

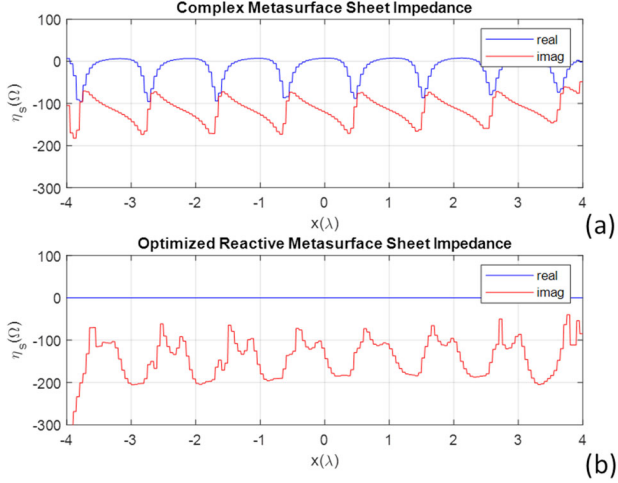

Figure 10. Metasurface layer 1 impedance sheet. (a) The complex-valued sheet resulting from the direct solution of the matrix equation. (b) The optimized purely reactive sheet.

If one takes the ratio of the beam peaks at $k_x/k_0=0$ (incident field) and $k_x/k_0=-0.94$ (reflected field) in Fig. 8, the result is 1.7 in agreement with (20). Note that (20) is a plane wave result, however the metasurface in Fig. 4 is finite. As a result, the scattered field spectrum of the finite metasurface will contain sidelobes at $k_x/k_0=0$, and the incident field will contain sidelobes at $k_x/k_0=-0.94$. By introducing surface waves (evanescent spectrum), we have ensured local power density conservation, resulting in a passive and lossless metasurface. Also note that the optimizer synthesized the power conserving scattered field amplitude despite the fact that the initial point had $E_r = E_0$.

The far field pattern of the optimized reactive metasurface is shown in Fig. 6 as the curve labeled 'MoM Opt React'. According to (19) and Fig. 7, the RMS difference between the patterns 'MoM Complex' and 'MoM Opt React' in Fig. 6 is 0.022 dB. The metasurface consisting of the optimized reactive $[\eta_{s1}]$ was modelled using the Finite Element Method (FEM) solver in COMSOL Multiphysics. The result is shown in Fig. 6 as the curve labeled 'COMSOL Opt React'. The full-wave verification is excellent.

Finally, the real part of the complex scattered near field was calculated and is plotted in Fig. 9. The scattered near field shows the reflected plane wave scattered to the desired angle. Note, the scattered field amplitude at the metasurface plane is 2.7 V/m and the incident plane wave's amplitude is 1V/m. When the complex-valued scattered field is added to the complex-valued incident field, the result is a reflected total field of amplitude 1.7 V/m as required for perfect anomalous reflection from (20). This is further evidenced by the plane wave in the region below the metasurface ($y < \lambda_0/23.62$) in Fig. 9. This scattered field has the same amplitude as the incident field but is 180 degrees out of phase. When added to the incident field, it produces a shadow zone or zero total field below the metasurface, as expected. Thus, the optimized purely reactive metasurface accomplishes perfect anomalous reflection.

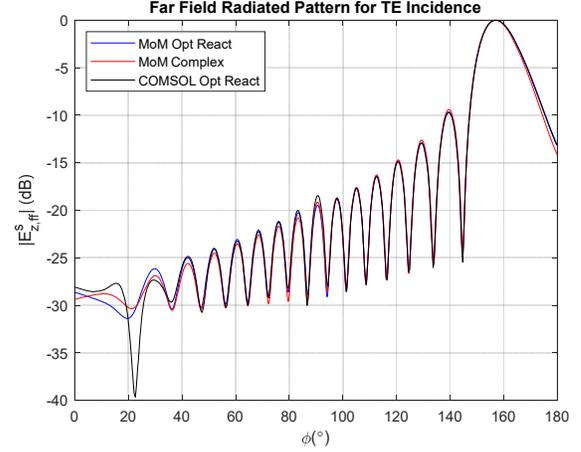

Figure 11. Scattered far field patterns. Shown is the far field pattern resulting from the complex-valued sheet calculated using the MoM (MoM Complex), the far field resulting from the optimized reactive sheet calculated using the MoM (MoM Opt React), and the far field pattern resulting from the optimized reactive sheet calculate in COMSOL Multiphysics (COMSOL Opt React).

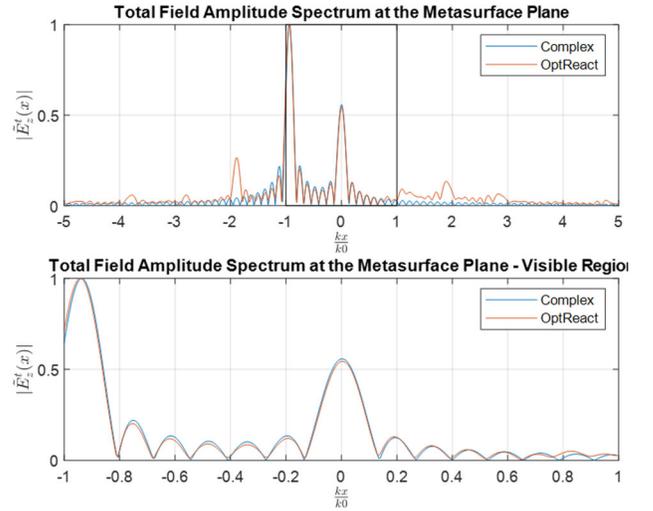

Figure 12. Total field amplitude spectra. (a) The visible and invisible spectrum. (b) Only the visible spectrum.

*Non-Uniqueness of Optimized Reactive Metasurface*

If instead the initial scattered field amplitude is chosen to satisfy, $E_r = E_0\sqrt{\cos(\theta_i)/\cos(\theta_r)}$, then the complex-valued metasurface sheet impedance will contain both positive and negative resistances such that the global integral of the power density along layer 1 results in a zero value [38]. As such, the metasurface then exhibits both loss and gain. This is evident in the impedance plots of Fig 10a where the resistance now has both positive and negative values. A comparison of Fig. 10a with Fig. 5a shows that both the real part and the imaginary part of the metasurface sheet impedance are very different, yet the far field pattern is the same (see Fig. 6 and Fig. 11). The fact that the scattered electric field amplitude follows (20) is evident in Fig. 12. Comparison with Fig. 8 shows that the complex-valued sheet already gets the scattered field amplitude for the reflected beam at $k_x/k_0$=-0.94 correct for local power conservation.

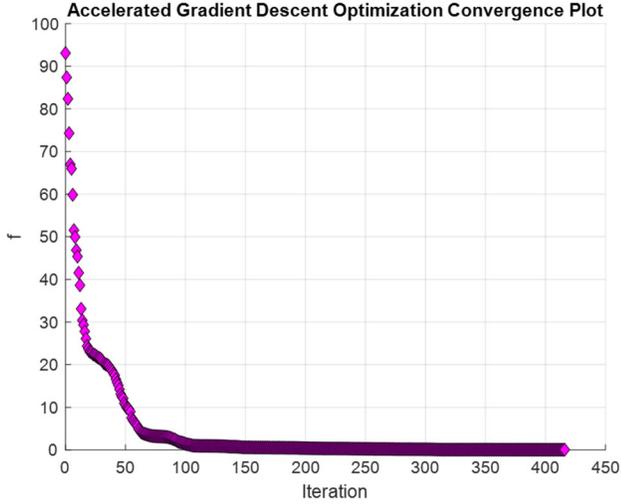

Figure 13. Optimization convergence plot.

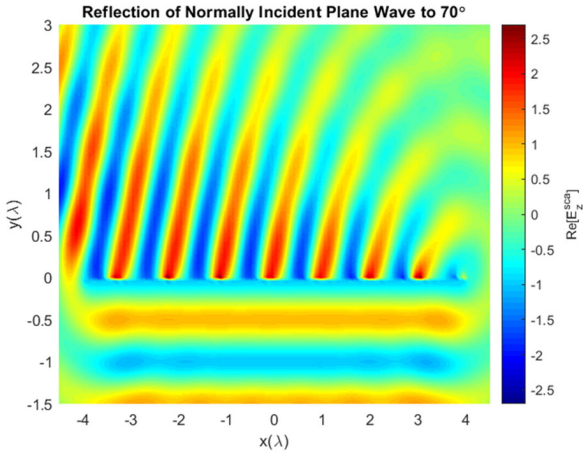

Figure 14. Real part of the complex scattered electric near field.

The same optimization technique was applied only this time the initial point is different. It is constructed from the imaginary part of Fig. 10a rather than the imaginary part of Fig. 5a. The optimization converged in 416 iterations to a value of $f$=0.068 (see Fig. 13) because the step size $\alpha$ reduced below $10^{-10}$ (see section II.E.2). This time, the optimization took 6 hours and 35 minutes running on the same 25 total cores of a machine with dual 3.0 GHz Intel Xeon Gold 6154 CPUs. The final optimized, purely reactive sheet is shown in Fig. 10b. It must be noted that this optimized reactive sheet is different than that in Fig. 5b, yet it leads to the same propagating scattered field within the visible spectrum in Fig. 12b. The invisible (evanescent) spectrum shown in Fig. 12a is different from that in Fig. 8a. This indicates that a different set of surface waves can lead to the same scattered far field (Fig. 6 and Fig. 11) and radiative near field (Fig. 9 and Fig. 14). Thus, the optimized reactive sheet is not unique. Different initial points can lead to different acceptable solutions with the same cost function $f$. This means that the cost function $f$ contains multiple local minima which all satisfy the convergence criteria in section II.E. Each initial set of reactances require different surface waves to achieve passivity. The key is finding the solution which requires the least amount of evanescent spectrum or lowest evanescent wavenumbers, as this solution will be less sensitive to loss and should exhibit wider bandwidth. With this in mind, a term can be added to the cost function $f$ defined in (19) which minimizes the necessary range of tangential wavenumbers in the invisible region.

*Patterning of Metallic Cladding and Full-Wave Simulation*

In this section, we pattern a metallic cladding to realize the optimized reactive sheet shown in Fig. 10b. This sheet is purely capacitive, has less evanescent content than the first design of Fig. 5b, and shows better agreement with the COMSOL full-wave simulations. The sheet impedances of the patterned geometries shown in Fig. 15 were extracted in a periodic environment using CST Microwave Studio, where the elements are placed at the interface between half spaces of air and dielectric. These geometries were then used to pattern the metasurface. The elements were placed in a parallel plate waveguide to emulate a geometry that is invariant in the z-direction. Full-wave simulations were carried out in COMSOL Multiphysics. The results are shown in Fig. 16a and Fig. 17. Although Fig. 17 shows the beam is reflected into the correct direction, Fig. 16a shows that the metasurface does not produce the same far field pattern. The incident plane wave specularly reflects, as can be seen by the beam peak at $\phi$=90°. This is also visible in Fig. 17 where the reflected plane wave in the normal direction causes perturbations to the wavefronts (these are not evident in Fig. 14). This is most likely because the extractions of the sheet impedances of the patterned geometries in Fig. 15 were performed in a locally periodic environment, whereas the optimized reactive sheet in Fig. 10b varies non-adiabatically in some places. The patterned metallic cladding was then finely tuned by optimizing the geometrical parameters $t$ (shown in Fig. 15) of all elements in the cladding simultaneously in COMSOL, in order to obtain closer agreement with the homogenized sheet version using a gradient descent optimization. Here, the Woodbury matrix identity cannot be used since the function $f$ in (19) is evaluated in COMSOL rather than by the method of moments. Note, in (19), $|E_{farfield}(\phi)|_{complex\ sheet}$ is replaced by $|E_{farfield}(\phi)|_{COMSOL\ patterned}$. A direct/brute force optimization of the patterned metallic claddings is not possible as each component of the gradient in (7) requires 7 minutes of computation time. This would result in a computation time of 18 hours and 40 minutes for full 160 component gradient of one optimization iteration. In comparison, using our algorithm, the full gradient is calculated in approximately 8 seconds since a fast matrix inversion rather than a volumetric finite element solution of finely patterned metallic cladding is required (also note the COMSOL simulations cannot be ran in parallel as the moment method solutions can). This shows the benefit of the homogenization and of our accelerated optimization method used in the design of the metasurface. Nonetheless, since we begin with a near-optimal solution, few iterations are required, and thus direct fine tuning of the patterned metallic cladding is made possible.

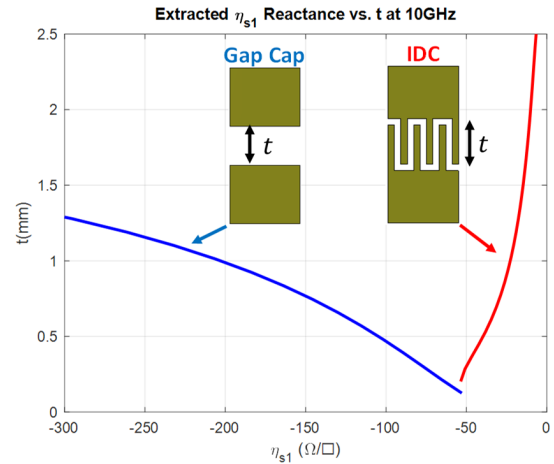

Figure 15. Patterned metallic cladding geometry sheet reactance versus geometrical parameter. Shown is the reactance versus geometrical parameter for both the Gap Capacitor (Gap Cap) and the Interdigitated Capacitor (IDC).

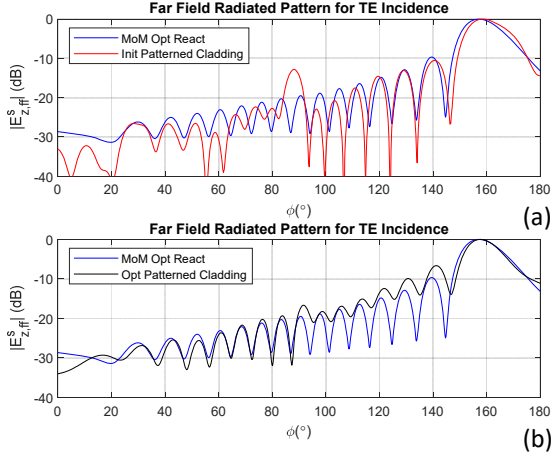

Figure 16. COMSOL Multiphysics full-wave simulations results of the metasurface made from patterned metallic claddings. Shown is the far field pattern of the optimized reactive sheet design calculated using the method of moments (MoM Opt React), *(a)* the far field pattern of the original cladding associated with the optimized reactive sheet design and calculated from COMSOL Multiphysics (COMSOL Cladding), and *(b)* the far field pattern of the optimized cladding in COMSOL Multiphysics (COMSOL Opt Cladding).

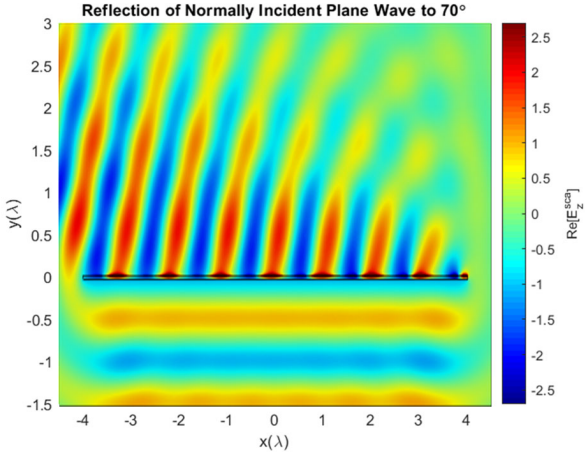

Figure 17. COMSOL Multiphysics full-wave simulations results of the metasurface made from patterned metallic claddings. Shown is the real part of the near field.

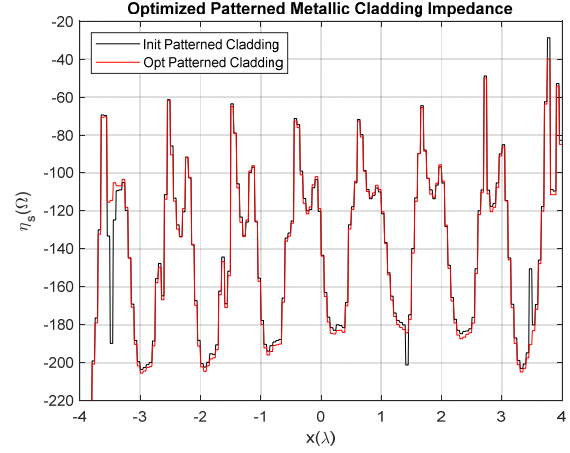

Figure 18. Patterned Metallic Cladding Optimization Result.

The results of the optimization are shown in Fig. 16b. The optimization completed in 16 days. The scattered far field shows closer agreement with that of the optimized reactive sheet impedance. The optimized cladding is compared to the initial cladding in Fig. 18. The plot was obtained using Fig. 15 but in the reverse way taking the geometrical parameter $t$ for each of the claddings and determining the equivalent periodic environment sheet impedance. This allows one to plot the change to the patterned cladding due to the fine tuning optimization. The root-mean-square (RMS) change in the impedances of the elements of the optimized cladding is 7.5Ω, indicating that impedances were off by on average 7.5Ω due to the assumption of local periodicity used during the extraction process.

V. CONCLUSION

A fast optimization method for finite non-periodic metasurfaces, containing hundreds of unknowns, modeled using integral equations is presented. The fast optimization method is based on a Woodbury Matrix Identity acceleration technique for calculating the gradient vector in gradient descent optimization methods. The technique leads to a 26.5 times improvement in computation time per iteration. The presented optimization algorithm has been used to design a wide-angle reflecting, lossless metasurface. The optimization results in different local power conserving design when initialized at different initial points. Thus, the optimized, purely reactive designs are not unique. This indicates that the required surface waves needed to obtain local power conservation depend on the initial reactances fed to the optimizer.